\documentclass[12pt]{article}
\usepackage{amsmath}
\usepackage{amscd}
\usepackage{mathrsfs}
\usepackage{amsfonts}
\usepackage{amsmath,amssymb}
\baselineskip 5.5pt \lineskip 5.5pt \numberwithin{equation}{section}

\def \Z{\hbox{$Z\hskip -5.2pt Z$}}

\def \C{\hbox{$C\hskip -5pt \vrule height 6pt depth 0pt \hskip 6pt$}}

\def\qed{\ \ \ifhmode\unskip\nobreak\fi\ifmmode\ifinner
         \else\hskip5pt\fi\fi
 \hbox{\hskip5pt\vrule width4pt height6pt depth1.5pt\hskip 1 pt}}
\def\a{\alpha}

\def\b{\beta}
\def\d{\delta}

\def\cl{\centerline}

\def\vs{\vspace*}

\def\C{\mathbb{C}}

\def\Z{\mathbb{Z}}

\newtheorem{theo}{Theorem}[section]
\newtheorem{lemm}[theo]{Lemma}

\newtheorem{defi}[theo]{Definition}
\newtheorem{prop}[theo]{Proposition}

\begin{document}

\begin{center}
{\Large \textbf{Intermediate series module of the Two Parameters Deformed Virasoro Algebra}
\noindent\footnote{Supported by the National Science Foundation of
China (Nos. 11047030 and 11771122).
}} \vs{6pt}
\end{center}

\cl{Wen Zhou, Yongsheng Cheng}

\cl{ \small School
of Mathematics and Statistics, Henan
University, Kaifeng 475004, China} \vs{6pt}

\vs{6pt}

\noindent{\it \textbf{Abstract.}~}
In this paper, we construct a class of Harish-Chandra modules of the two parameters deformed Virasoro algebra and classify indecomposanle Harish-Chandra module of an intermediate series.\\
\noindent{\it \textbf{Keywords:}~}
the two parameters deformed Virasoro algebra; intermediate series module; indecomposanle Harish-Chandra module.

{\small
\parskip .005 truein
\baselineskip 10pt \lineskip 10pt
\cl{\bf\S1. \ Introduction}
  Recently, one of the most modern trends in mathematics that how to classify all Harish-Chandra modules, especially some Lie algebras related to the Virasoro algebra.\par
  In \cite{[IK]}, Kaplansky studied a class of Harish-Chandra modules of the centerless Virasoro algebra in 1982, and classified all indecomposable Harish-Chandra modules of the Virasoro algebra in 1985(see \cite{[IL]}). K. Liu gave a class of Harish-Chandra modules of Virasoro algebra (centerless) of Hom type (q is not a unit root), and provided the classification of all its indecomposable Harish-Chandra modules in 1995(see \cite{[KL]}). Z. Zhang, G. Zhang and Y. Jia constructed a class of Harish-Chandra modules with multiplicity $\leq$ 1 of the two parameters deformation a Virasoro algebra and proved a classification theorem(see \cite{[ZZ]}). In \cite{[LS]}, J. Li and Y. Su proved that an irreducible weight module with finite-dimensional weight spaces over the Schr\"{o}dinger-Virasoro algebras is a highest/lowest weight module or a uniformly bounded module. Furthermore, indecomposable modules of the intermediate series over these algebras are completely determined. K. Zhao and R. Lv classified all of the Harish-Chandra modules on the twisted Heisenberg-Virasoro algebra, but the caculation is too complicated and the computational content is too large(see \cite{[LZ]}). In \cite{[DL]}, D. Liu provided a uniform method to thoroughly classify Harish-Chandra modules over some Lie algebras related to the Virasoro algebras.\par
  This paper is organized as follows. In Section 2, we recall some basic definitions of $V_{p,q}$. In Section 3,
   we provide a class of Harish-Chandra modules of the two parameters deformed Virasoro algebra, and proves a classsification theorem. Moreover, we classify indecomposable Harish-Chandra module of the two parameters deformed Virasoro algebra. In particular, after calculating the classification theorem, we found that the quantum integer in this paper are consistent with it in the \cite{[ZZ]} after replacement.\\

\cl{\bf\S2.\ Notations and Preliminaries}
\setcounter{section}{2}\setcounter{theo}{0}\setcounter{equation}{0}
In \cite{[ELMS]}, E. Olivier et. al introduced the definition of Hom-Lie algebra $V_{p,q}$.
\begin{theo}\label{1}$^{\cite{[ELMS]}}$
The two parameters deformed Virasoro algebra $V_{p,q}=(\hat{L},\hat{\alpha})$, where
$\hat{L}$ has basis $\{L_n,C|n\in\Z\}$ and bracket relations:
\begin{align*}
[L_n,L_m]:&=(\frac{[n]}{p^{n}}-\frac{[m]}{p^{m}})L_{n+m}
+\d_{m+n,0}\frac{(q/p)^{-n}}{6(1+(q/p)^{n})}
\frac{[n-1]}{p^{n-1}}
\frac{[n]}{p^{n}}\frac{[n+1]}{p^{n+1}}C,\\
[\hat{L},C]:&=0,
\end{align*}
and $\hat{\alpha}:\hat{L}\longrightarrow \hat{L}$ is the endomorphism of $\hat{L}$ defined by\\
$$\hat{\alpha}(L_n)=((1+(q/p)^{n}))L_n,\ \hat{\alpha}(C)=C.$$
\end{theo}
Let us recall the definition of the module.
\begin{defi}
A vector space M is called a $V_{p,q}$-module if there exists a map:
\begin{align*}
V_{p,q}\times M &\rightarrow M\\
(X,v)&\mapsto Xv\ (\forall X\in V_{p,q}, v\in M)
\end{align*}
such that
\begin{align*}
(p^{-n}q^{n}L_{n}L_{m}-p^{-m}q^{m}L_{m}L_{n})v=\big(\frac{[m]_{p,q}}{p^{m}}-\frac{[n]_{p,q}}{p^{n}}
\big)L_{m+n}v,\ (\forall m,n\in \mathbb{Z})
\end{align*}
 where the quantum integer
\begin{align*}
&[n]_{p,q}=\frac{p^{n}-q^{n}}{p-q}.
\end{align*}
The $V_{p,q}$-module M is called a Harish-Chandra module if $M=\oplus_{\lambda\in \mathbb{C}}M_{\lambda}$ where $M_{\lambda}=\{v\in M|L_{0}v=\lambda v\}$.\\
\end{defi}

\begin{theo}\label{1}$^{\cite{[ZZ]}}$
$\forall (a,b)\in \mathbb{C}\times\mathbb{C}$, we construct a $V_{p,q}$-module $M_{a,b}$ as follows:
\begin{align}
M_{a,b}&=\oplus_{k\in \mathbb{Z}}\mathbb{C}v_{k},\notag\\
L_{n}(v_{k})&=(p^{-k}[k]_{p,q}-ap^{-k}q^{k}-bp^{-k-n}q^{k}[n]_{p,q})v_{k+n}. \  (\forall n,k\in \mathbb{Z})
\end{align}
A direct calculation can check that $M_{a,b}$ is a $V_{p,q}$-module. It is clear that $M_{a,b}^{\lambda}=\oplus_{\lambda\in \mathbb{C}}M_{a,b}^{\lambda}$, where $M_{a,b}^{\lambda}=\{v\in M_{a,b}|L_{0}v=\lambda v\}$. $\lambda$ is called a weight if $dim M_{a,b}^{\lambda}>0$, and $dim M_{a,b}^{\lambda}$ is called the multiplicity of $\lambda$ denoted by $mult\lambda$. We denote the weight set of $M_{a,b}$ by $P(M_{a,b})$. It is clear that $P(M_{a,b})=\{p^{-k}[k]_{p,q}-ap^{-k}q^{k}|k\in \mathbb{Z}\}$, and $mult(p^{-k}[k]_{p,q}-ap^{-k}q^{k})=1$ for all $k\in \mathbb{Z}$ if $a\neq -\frac{1}{p-q}$.
\end{theo}

\cl{\bf\S3.\ Intermediate series modules of the $V_{p,q}$ }
\setcounter{section}{3}\setcounter{theo}{0}\setcounter{equation}{0}

In this section, we give some basic propositions about $M_{a,b}$, then we study the intermediate series modules of the two parameters deformed Virasoro algebra $V_{p,q}$.

\begin{prop}
Let $a\neq -\frac{1}{p-q}$, then $M_{a,b}\simeq M_{a^{'},b^{'}}$ if and only if $a=-p^{-m}[m]_{p,q}+p^{-m}q^{m}a^{'}, b=p^{-m}q^{m}b^{'}$ for some $m\in \mathbb{Z}$.
\end{prop}
\noindent{\it \textbf{Proof.}~}
Set $M_{a,b}=\oplus_{i\in \mathbb{Z}}\mathbb{C}v_{i}$, $M_{a^{'},b^{'}}=\oplus_{i\in \mathbb{Z}}\mathbb{C}v_{i}^{'}$.\\
Let $\phi$ be the isomorphism from $M_{a,b}$ to $M_{a^{'},b^{'}}$. It follows from $P(M_{a,b})=P(M_{a^{'},b^{'}})$ that :
\begin{align*}
p^{-k}[k]_{p,q}-ap^{-k}q^{k}=p^{-k^{'}}[k^{'}]_{p,q}-a^{'}p^{-k^{'}}q^{k^{'}},
\end{align*}
then
\begin{align*}
a=p^{k-k^{'}}\frac{q^{k-k^{'}}}{p^{k^{'}-k}}+p^{k-k^{'}}q^{k^{'}-k}a^{'},
\end{align*}
let $k^{'}-k=m$, we get\\
\begin{align*} a=-p^{-m}[m]_{p,q}+p^{-m}q^{m}a^{'}
\end{align*}
for some $m\in \mathbb{Z}$. By $\phi(L_{n}v_{k})=L_{n}(\phi(v_{k}))$, we see that: $b=p^{-m}q^{m}b^{'}$.\\
Define a map $\phi:M_{a,b}\rightarrow M_{a^{'},b^{'}}$ by $v_{i}\mapsto v_{i+m}^{'}$, using $\phi(L_{n}v_{k})=L_{n}(v_{k+m}^{'})$, it is clear that $\phi$ is a module isomorphism.
\hfill$\Box$\vskip7pt

\begin{prop}
Let $a\neq -\frac{1}{p-q}$, $M_{a,b}=\oplus_{i\in \mathbb{Z}}\mathbb{C}v_{i}$ is reducible if and only if $a=-p^{-m}[m]_{p,q}, b=-p^{-m}q^{m} or$ 0 for some $m\in \mathbb{Z}$.\\
\noindent{\it \textbf{Proof.}~}
Suppose that $U$ is a nonzero proper submodule of $M_{a,b}$. By the action of $L_{0}$, we can find that there exists a $k\in \mathbb{Z}$ such that $v_{k}\in U$. Hence, one can find that there exist some $n\in \mathbb{Z}$ such that $p^{-k}[k]_{p,q}-ap^{-k}q^{k}-bp^{-k-n}q^{k}[n]_{p,q}=0$. We have the following two cases:\\
\textbf{Case 1.}
Then $n$ is unique. Then $v_{i}\in \mathbb{C}L_{i-k}v_{k}\subset U$ for all $i\neq n+k$. Let $m=-n-k$, we see that $U=\oplus_{i\neq -m}\mathbb{C}v_{i}$. Using $L_{i-m}v_{-i}=0\ (\forall i\neq m)$, we get
\begin{align*}
p^{i}[-i]_{p,q}-ap^{i}q^{-i}-bp^{m}q^{-i}[i-m]_{p,q}=0,
\end{align*}
then: $a=-p^{-m}[m]_{p,q}$, $b=-p^{-m}q^{-m}$.\\
\textbf{Case 2.}
The $n$ is not unique. Suppose that there exist $n_{1}, n_{2}\in \mathbb{Z}$ such that
\begin{align*}
p^{-k}[k]_{p,q}-ap^{-k}q^{k}-bp^{-k-n_{1}}q^{k}[n_{1}]_{p,q}=0,\\
p^{-k}[k]_{p,q}-ap^{-k}q^{k}-bp^{-k-n_{2}}q^{k}[n_{2}]_{p,q}=0,
\end{align*}
thus, $b=0, a=q^{-k}[k]_{p,q}$, and $U=\mathbb{C}v_{k}$.\\
Suppose that $a=-p^{-m}[m]_{p,q}, b=-p^{-m}q^{m} or$ 0, it is a simple argument.
\hfill$\Box$\vskip7pt
\end{prop}
Before giving the classification theorem we need the following lemma.
\begin{lemm}
$V$ is generated by $L_{\pm1}, L_{\pm2}$.\\
\end{lemm}

\begin{theo}
Let $M=\oplus_{k\in \mathbb{Z}}\mathbb{C}v_{k}$ be an irreducible Harish-Chandra module over $V_{p,q}$ such that $L_{n}v_{k}\in \mathbb{C}v_{n+k},\  Ann_{M}(L_{\pm1})=0$ and $-\frac{1}{p-q}\bar{\in}P(M)$. Then there exist $a,b\in \mathbb{C}$, such that $M\simeq M_{a,b}$.
\end{theo}
\noindent{\it \textbf{Proof.}~}
It is clear that $a\neq -\frac{1}{p-q}$. Since $v_{k}\in \mathbb{C}L_{1}^{k}v_{0}$ for all $k>0$, and $v_{k}\in \mathbb{C}L_{-1}^{-k}v_{0}$ for all $k<0$, we see that
\begin{align}
L_{0}v_{k}=(p^{-k}[k]_{p,q}-ap^{-k}q^{k})v_{k},\label{3.1}\ k\in \mathbb{Z}.
\end{align}
Set $L_{1}v_{i}=a_{i}v_{i+1},\ L_{-1}v_{i}=b_{i}v_{i-1},\ (i\in \mathbb{Z},\  a_{i},b_{i}\in \mathbb{C}^{*})$. Putting $n=1, m=-1$ in (\ref{3.1}), we see that
\begin{align}
(p^{-1}qL_{1}L_{-1}-pq^{-1}L_{-1}L_{1})v_{i}=pq[-2]_{p,q}L_{0}v_{i},\label{3.2}
\end{align}
thus
\begin{align*}
p^{-1}qa_{i-1}b_{i}-pq^{-1}a_{i}b_{i+1}=pq[-2]_{p,q}(p^{-i}[i]_{p,q}-ap^{-i}q^{i}),\ i\in \mathbb{Z}.
\end{align*}
Let $b$ be one of the solutions of $a_{0}b_{1}=(-a-xp^{-1}[1]_{p,q})(p^{-1}[1]_{p,q}-ap^{-1}q-xq[-1]_{p,q})$, then by (\ref{3.3}), we get
\begin{align}
b_{i+1}a_{i}=&(p^{-i-1}[i+1]_{p,q}-ap^{-i-1}q^{i+1}-bp^{-i}q^{i+1}[-1]_{p,q})\notag\\
&(p^{-i}[i]_{p,q}-ap^{-i}q^{i}-bp^{-i-1}q^{i}[1]_{p,q}),\label{3.3}\ (\forall i\in \mathbb{Z}).
\end{align}
We fix $v_{0}$ and determine $v$'s by
\begin{align}
L_{1}v_{i}=(p^{-i}[i]_{p,q}-ap^{-i}q^{i}-bp^{-i-1}q^{i}[1]_{p,q})v_{i+1},\label{3.4}\ (\forall i\in \mathbb{Z}).
\end{align}
From (\ref{3.4}) and (\ref{3.5}) we decuce
\begin{align}
L_{-1}v_{i}=(p^{-i}[i]_{p,q}-ap^{-i}q^{i}-bp^{-i+1}q^{i}[-1]_{p,q})v_{i-1},\label{3.5}\ (\forall i\in \mathbb{Z}).
\end{align}
Set
\begin{align}
L_{2}v_{j}=(f(j)+p^{-j}[j]_{p,q}-ap^{-j}q^{j}-bp^{-j-2}q^{j}[2]_{p,q})v_{j+2},\label{3.6}\\ L_{-2}v_{j}=(g(j)+p^{-j}[j]_{p,q}-ap^{-j}q^{j}-bp^{-j+2}q^{j}[-2]_{p,q})v_{j-2},\label{3.7}
\end{align}
where $f(j),g(j)\in \mathbb{C}$.
Using
\begin{align*}
(p^{-2}q^{2}L_{2}L_{-1}-pq^{-1}L_{-1}L_{2})v_{j}&=(\frac{[-1]_{p,q}}{p^{-1}}-\frac{[2]_{p,q}}{p^{2}})L_{1}v_{j},\\
(p^{2}q^{-2}L_{-2}L_{1}-p^{-1}qL_{1}L_{-2})v_{j}&=(\frac{[1]_{p,q}}{p}-\frac{[-2]_{p,q}}{p^{-2}})L_{-1}v_{j},
\end{align*}
we deduce
\begin{align}
q^{-3}f(j)(p^{-j-2}[j+2]_{p,q}-ap^{-j-2}q^{j+2}-bp^{-j-1}q^{j+2}[-1]_{p,q})\notag\\
=p^{-3}f(j-1)(p^{-j}[j]_{p,q}-ap^{-j}q^{j}-bp^{-j+1}q^{j}[-1]_{p,q}),\label{3.8}
\end{align}
and
\begin{align}
p^{3}g(j+1)(p^{-j}[j]_{p,q}-ap^{-j}q^{j}-bp^{-j-1}q^{j}[1]_{p,q})\notag\\
=q^{3}g(j)(p^{-j+2}[j-2]_{p,q}-ap^{-j+2}q^{j-2}-bp^{-j+1}q^{j-2}[1]_{p,q}).\label{3.9}
\end{align}
Then
\begin{align}
f(j)&=\frac{p^{-3j}q^{3j}F}{(p^{-j-2}[j+2]_{p,q}-ap^{-j-2}q^{j+2}-bp^{-j-1}q^{j+2}[-1]_{p,q})
(p^{-j-1}[j+1]_{p,q}-ap^{-j-1}q^{j+1}-bp^{-j}q^{j+1}[-1]_{p,q})},\label{3.10}\\
g(j)&=\frac{p^{-3j}q^{3j}G}{(p^{-j+2}[j-2]_{p,q}-ap^{-j+2}q^{j-2}-bp^{-j+1}q^{j-2}[1]_{p,q})
(p^{-j+1}[j-1]_{p,q}-ap^{-j+1}q^{j-1}-bp^{-j}q^{j-1}[1]_{p,q})},\label{3.11}
\end{align}
where
\begin{align*}
F&=f(0)(p^{-2}[2]_{p,q}-ap^{-2}q^{2}-bp^{-1}q^{2}[-1]_{p,q})
(p^{-1}-ap^{-1}q-bq[-1]_{p,q}),\\
G&=g(0)(p^{2}[-2]_{p,q}-ap^{2}q^{-2}-bpq^{-2})
(p[-1]_{p,q}-apq^{-1}-bq^{-1}).
\end{align*}
Let $x=q^{-j}[j]_{p,q}$, then we have
\begin{equation}
\begin{cases}
p^{j}q^{-j}f(j)=\frac{F}{(x+p^{-2}[2]_{p,q}-ap^{-2}q^{2}-bp^{-1}q^{2}[-1]_{p,q})
(x+p^{-1}-ap^{-1}q^{1}-bq[-1]_{p,q})},\\
p^{j}q^{-j}f(j-2)=\frac{p^{6}q^{-6}F}{(x-a-bp[-1]_{p,q})
(x+p[-1]_{p,q}-apq^{-1}-bp^{2}q[-1]_{p,q})},\label{3.12}\\
p^{j}q^{-j}g(j)=\frac{G}{(x+p^{2}[-2]_{p,q}-ap^{2}q^{-2}-bpq^{-2})
(x+p[-1]_{p,q}-apq-bq^{-1})},\\
p^{j}q^{-j}g(j+2)=\frac{p^{-6}q^{6}G}{(x-a-bp^{-1})
(x+p^{-1}-ap^{-1}q-bp^{-2}q)},
\end{cases}
\end{equation}
By $(p^{-2}q^{2}L_{2}L_{-2}-p^{2}q^{-2}L_{-2}L_{2})v_{j}=p^{-2}q^{2}[-4]_{p,q}L_{0}v_{j}$, we get
\begin{align}
&p^{-2}q^{2}\big((p^{-j}[j]_{p,q}-ap^{-j}q^{j}-bp^{-j+2}q^{j}[-2]_{p,q})f(j-2)+\notag\\
&(p^{-j+2}[j-2]_{p,q}-ap^{-j+2}q^{j-2}-bp^{-j}q^{j-2}[2]_{p,q})g(j)+f(j-2)g(j)\big)=\label{3.13}\\
&p^{2}q^{-2}\big((p^{-j-2}[j+2]_{p,q}-ap^{-j-2}q^{j+2}-bp^{-j}q^{j+2}[-2]_{p,q})f(j)+\notag\\
&(p^{-j}[j]_{p,q}-ap^{-j}q^{j}-bp^{-j-2}q^{j}[2]_{p,q})g(j+2)+f(j)g(j+2)\big).\notag
\end{align}
Using (\ref{3.13}) and (\ref{3.14}), we have
\begin{align}
&p^{-4}q^{4}(x+p^{-2}[2]_{p,q}-ap^{-2}q^{2}-bp^{-1}q^{2}[-1]_{p,q})
(x+p^{-1}-ap^{-1}q^{1}-bq[-1]_{p,q})\notag\\
&(x-a-bp^{-1})
(x+p^{-1}-ap^{-1}q-bp^{-2}q)\big(p^{6}q^{-6}F(x+p^{2}[-2]_{p,q}-ap^{2}q^{-2}-bpq^{-2})\notag\\
&(x+p[-1]_{p,q}-apq-bq^{-1})(x-a-bp^{2}[-2]_{p,q})+G(x-a-bp[-1]_{p,q})\notag\\
&(x+p[-1]_{p,q}-apq^{-1}-bp^{2}q[-1]_{p,q})(x+p^{2}[-2]_{p,q}-ap^{2}q^{-2}-bpq^{-2}[2]_{p,q})+p^{6}q^{-6}FG\big)\notag\\
=&(x-a-bp[-1]_{p,q})
(x+p[-1]_{p,q}-apq^{-1}-bp^{2}q[-1]_{p,q})(x+p^{2}[-2]_{p,q}-ap^{2}q^{-2}-bpq^{-2})\label{3.14}\\
&(x+p[-1]_{p,q}-apq^{-1}-bq^{-1})\big(F(x-a-bp^{-1})
(x+p^{-1}-ap^{-1}q-bp^{-2}q)\notag\\
&(x+p^{-2}[2]_{p,q}-ap^{-2}q^{2}-bq^{2}[-2]_{p,q})+p^{-6}q^{6}G(x+p^{-2}[2]_{p,q}-ap^{-2}q^{2}-bp^{-1}q^{2}[-1]_{p,q})\notag\\
&(x+p^{-1}-ap^{-1}q-bq[-1]_{p,q})(x-a-bp^{-2}[2]_{p,q})+p^{-6}q^{6}FG
\big)\notag
\end{align}
It is clear that (\ref{3.15}) is a plooynomial isentity. Comparing the cofficients of terms of the highest degree of the two sides, we obtain that $F=-p^{-6}q^{6}G$, if $F=0$, we see that $M\simeq M_{a,b}$ are required; if $F\neq 0$, we have
\begin{equation*}
\begin{cases}
f_{1}(x)=&x-a-bp^{-1},\\
f_{2}(x)=&x+p^{-1}-ap^{-1}q-bp^{-2}q,\\
f_{3}(x)=&x+p^{-2}[2]_{p,q}-ap^{-2}q^{2}-bp^{-1}q^{2}[-1]_{p,q},\\
f_{4}(x)=&x+p^{-1}-ap^{-1}q^{1}-bq[-1]_{p,q},\\
f_{5}(x)=&(x-a-bp^{-1})
(x+p^{-1}-ap^{-1}q-bp^{-2}q)(x-a-bp^{2}[-2]_{p,q})-\\
&(x-a-bp[-1]_{p,q})(x+p[-1]_{p,q}-apq^{-1}-bp^{2}q[-1]_{p,q})\\
&(x+p^{2}[-2]_{p,q}-ap^{2}q^{-2}-bpq^{-2}[2]_{p,q})-p^{6}q^{-6}F,
\end{cases}
\end{equation*}
and\\
\begin{equation*}
\begin{cases}
g_{1}(x)=&x-a-bp[-1]_{p,q},\\
g_{2}(x)=&x+p[-1]_{p,q}-apq^{-1}-bp^{2}q[-1]_{p,q},\\
g_{3}(x)=&x+p[-1]_{p,q}-apq^{-1}-bq^{-1},\\
g_{4}(x)=&x+p^{2}[-2]_{p,q}-ap^{2}q^{-2}-bpq^{-2},\\
g_{5}(x)=&(x+p^{2}[-2]_{p,q}-ap^{2}q^{-2}-bpq^{-2})
(x+p[-1]_{p,q}-apq-bq^{-1})\\
&(x+p^{-2}[2]_{p,q}-ap^{-2}q^{2}-bq^{2}[-2]_{p,q})-
(x+p^{-2}[2]_{p,q}-ap^{-2}q^{2}-bp^{-1}q^{2}[-1]_{p,q})\\
&(x+p^{-1}-ap^{-1}q^{1}-bq[-1]_{p,q})
(x-a-bp^{-2}[2]_{p,q})-F.
\end{cases}
\end{equation*}
Thus, (\ref{3.15}) turns into
\begin{align}
f_{1}(x)f_{2}(x)f_{3}(x)f_{4}(x)f_{5}(x)=g_{1}(x)g_{2}(x)g_{3}(x)g_{4}(x)g_{5}(x).\label{3.15}
\end{align}
Through direct calculation we can find that
\begin{align}
deg f_{5}(x)\leq 1, deg g_{5}(x)\leq 1.\label{3.16}
\end{align}
About $f_{i}(x), g_{i}(x),(1\leq i\leq 4)$, we get:
\begin{align}
&f_{1}(x)=g_{1}(x)\Leftrightarrow b=0,\notag\\
&f_{2}(x)=g_{1}(x)\Leftrightarrow a-b\frac{q^{-1}+p^{-2}q}{1-p^{-1}q}=-\frac{1}{p-q},\notag\\
&f_{1}(x)=g_{2}(x)\Leftrightarrow a-b\frac{pq^{-2}-q^{-1}}{pq^{-1}-1}=-\frac{1}{p-q},\notag\\
&f_{2}(x)=g_{2}(x)\Leftrightarrow a-b\frac{p^{-1}q+pq^{-2}}{pq^{-1}-p^{-1}q}=-\frac{1}{p-q},\notag\\
&f_{1}(x)=g_{3}(x)\Leftrightarrow a-b\frac{p^{-1}-q^{-1}}{pq^{-1}-1}=-\frac{1}{p-q},\notag\\
&f_{2}(x)=g_{3}(x)\Leftrightarrow a-b\frac{p^{-2}q-q^{-1}}{pq^{-1}-p^{-1}q}=-\frac{1}{p-q},\notag\\
&f_{1}(x)=g_{4}(x)\Leftrightarrow a-bp^{-1}=-\frac{1}{p-q},\notag\\
&f_{2}(x)=g_{4}(x)\Leftrightarrow a-bp^{-1}=-\frac{1}{p-q},\notag\\
&f_{3}(x)=g_{1}(x)\Leftrightarrow a-b\frac{q^{-1}-p^{-2}q}{1-p^{-2}q^{2}}=-\frac{1}{p-q}, \label{3.17} \\
&f_{4}(x)=g_{1}(x)\Leftrightarrow a-b\frac{q^{-1}-p^{-1}}{1-p^{-1}q}=-\frac{1}{p-q},\notag\\
&f_{3}(x)=g_{2}(x)\Leftrightarrow a-b\frac{pq^{-2}-p^{-2}q}{pq^{-1}-p^{-2}q^{2}}=-\frac{1}{p-q},\notag\\
&f_{4}(x)=g_{2}(x)\Leftrightarrow a-b\frac{pq^{-2}-p^{-1}}{pq^{-1}-p^{-1}q}=-\frac{1}{p-q},\notag\\
&f_{3}(x)=g_{3}(x)\Leftrightarrow a-b\frac{-p^{-2}q-q^{-1}}{pq^{-1}-p^{-2}q^{2}}=-\frac{1}{p-q},\notag\\
&f_{4}(x)=g_{3}(x)\Leftrightarrow a-b\frac{-1}{p-q}=-\frac{1}{p-q},\notag\\
&f_{3}(x)=g_{4}(x)\Leftrightarrow a-b\frac{-pq^{-2}-p^{-2}q}{p^{2}q^{-2}-p^{-2}q^{2}}=-\frac{1}{p-q},\notag\\
&f_{4}(x)=g_{4}(x)\Leftrightarrow a-b\frac{-p^{-1}-pq^{-2}}{p^{2}q^{-2}-p^{-1}q}=-\frac{1}{p-q}.\notag
\end{align}
By $a\neq -\frac{1}{p-q}$ and (\ref{3.17}), one can find that $f_{i}(x),g_{i}(x)$ satisfy one of the following four cases:\\
(1) $f_{i}(x)\neq g_{i}(x),\ (\forall 1\leq i,j\leq 4)$;\\
(2) $f_{1}(x)= g_{1}(x), f_{i}(x)\neq g_{i}(x)$ for others;\\
(3) $f_{1}(x)= f_{2}(x)= g_{3}(x)= g_{4}(x),
f_{i}(x)\neq g_{i}(x)$ for others;\\
(4) $f_{3}(x)= f_{4}(x)= g_{1}(x)= g_{2}(x),
f_{i}(x)\neq g_{i}(x)$ for others.\\
If $f_{i}(x)\neq g_{i}(x)$, then $(f_{i}(x),g_{i}(x))=1$.
$g_{1}(x)g_{2}(x)g_{3}(x)g_{4}(x)$ divides $f_{5}(x)$, $f_{1}(x)f_{2}(x)f_{3}(x)f_{4}(x)$ divides $g_{5}(x)$. By (\ref{3.17}), we have
\begin{align}
 f_{5}(x)=g_{5}(x)=0.\label{3.18}
\end{align}
Thus
\begin{align}
F&=(x-a-bp^{-1})
(x+p^{-1}-ap^{-1}q-bp^{-2}q)(x+p^{-2}[2]_{p,q}-ap^{-2}q^{2}-bq^{2}[-2]_{p,q})\notag\\
&-(x+p^{-2}[2]_{p,q}-ap^{-2}q^{2}-bp^{-1}q^{2}[-1]_{p,q})
(x+p^{-1}-ap^{-1}q-bq[-1]_{p,q})(x-a-bp^{-2}[2]_{p,q})\label{3.19}
\end{align}
\begin{align}
G&=(x-a-bp[-1]_{p,q})
(x+p[-1]_{p,q}-apq^{-1}-bp^{2}q[-1]_{p,q})(x+p^{2}[-2]_{p,q}-ap^{2}q^{-2}-bpq^{-2}[2]_{p,q})\notag\\
&-(x+p^{2}[-2]_{p,q}-ap^{2}q^{-2}-bpq^{-2})
(x+p[-1]_{p,q}-apq^{-1}-bq^{-1})(x-a-bp^{-2}[2]_{p,q})\label{3.20}
\end{align}
Putting into (\ref{3.7}),(\ref{3.8}), using (\ref{3.11}) and (\ref{3.12}), we see that
\begin{align}
L_{2}(v_{j})=&\frac{(p^{-j}[j]_{p,q}-ap^{-j}q^{j}-bp^{-j-1}q^{j})(p^{-j-1}[j+1]_{p,q}-ap^{-j-1}q^{j+1}-bp^{-j-2}q^{j+1}))}
{(p^{-j-2}[j+2]_{p,q}-ap^{-j-2}q^{j+2}-bp^{-j-1}q^{j+2}[-1]_{p,q})}\notag\\
&\frac{(p^{-j-2}[j+2]_{p,q}-ap^{-j-2}q^{j+2}-bp^{-j}q^{j+2}[-2]_{p,q})}{(p^{-j-1}[j+1]_{p,q}-ap^{-j-1}q^{j+1}-bp^{-j}q^{j+1}[-1]_{p,q})},\label{3.21}\\
L_{-2}(v_{j})=&\frac{(p^{-j}[j]_{p,q}-ap^{-j}q^{j}-bp^{-j-1}q^{j}[-1]_{p,q})(p^{-j+1}[j-1]_{p,q}-ap^{-j+1}q^{j-1}-bp^{-j+2}q^{j-1}[-1]_{p,q}))}
{(p^{-j+2}[j-2]_{p,q}-ap^{-j+2}q^{j-2}-bp^{-j+1}q^{j-2})}\notag\\
&\frac{(p^{-j+2}[j-2]_{p,q}-ap^{-j+2}q^{j-2}-bp^{-j}q^{j-2}[-2])}{(p^{-j+1}[j-1]_{p,q}-ap^{-j+1}q^{j-1}-bp^{-j}q^{j-1})}\label{3.22}
\end{align}
We choose $h_{0}=1, h_{i}\in \mathbb{C}^{*}\ (\forall i\in \mathbb{Z})$ such that
\begin{align*}
\frac{h_{i}}{h_{i+1}}=\frac{p^{-i-1}[i+1]_{p,q}-ap^{-i-1}q^{i+1}-bp^{-i}q^{i+1}[-1]_{p,q}}
{p^{-i}[i]_{p,q}-ap^{-i}q^{i}-bp^{-i-1}q^{i}}
\end{align*}
Set $u_{j}=h_{j}v_{j}$, one can find
\begin{align*}
L_{n}u_{j}=(p^{-j}[j]_{p,q}-ap^{-j}q^{j}-b^{'}p^{-j-n}q^{j}[n]_{p,q})u_{n+j},
\end{align*}
where $j\in \mathbb{Z}, n=\pm1,\pm2$,\ and $b^{'}=1-a(p-q)-b$ is another solution of the equation $a_{0}b_{1}=(-a-xp^{-1}[1]_{p,q})(p^{-1}[1]_{p,q}-ap^{-1}q-xq[-1]_{p,q})$. \\ By Lemma we can prove $M\simeq M_{a,b^{'}}$.
\hfill$\Box$\vskip7pt

$U_{q}(sl_{2})$ is the simplest quantized universal enveloping algebra, suppose the elements $K,E,F,K^{-1}$ span $U_{q}(sl_{2})$. $T_{\omega l}$ is a representation space for it, define operators $T_{\omega l}(K)$, $T_{\omega l}(E)$, $T_{\omega l}(F)$ acting on $V_{l}$ by
\begin{align*}
&T_{\omega l}(K)e_{m}=\omega q^{2m}e_{m},\\
&T_{\omega l}(E)e_{m}=([l-m][l+m+1])^{\frac{1}{2}}e_{m+1},\\
&T_{\omega l}(F)e_{m}=\omega([l+m][l-m+1])^{\frac{1}{2}}e_{m-1},
\end{align*}
where $\omega\in\{-1,+1\}$, $[n]=\frac{q^{n}-q^{-n}}{q-q^{-1}}$, $V_{l}$ is a 2$l$+1 dimensional complex vector space with basis $e_{m},m=-l,-l+1,\cdots,l$. From the representation of the $U_{q}(sl_{2})$, we have
the characteristic roots of $K$ on the $e_{m}^{,}$ change by a factor of $\pm q^2$, and $T_{\omega l}(F)T_{\omega l}(E)e_{m}=\omega([l-m][l+m+1])e_{m}$, if we set $q^{-m}[m]=x$, then $[l-m][l+m+1]$ is a quadratic polynomial in $x$, so the coefficient of $e_{m}$ in $q^{-2m}T_{\omega l}(F)T_{\omega l}(E)(e_{m})$ is a quadratic polynomial in $x$. This fact will be used repeatly in the following proof.

 Then we can consider the case where $L_{-1}$ or $L_{1}$ annihilates. \\
 \textbf{Case I:}\\
 Using the inverted module of $V_{p,q}$, we can suppose that $L_{-1}$ annihilates some $v_{i}$. Let $L_{-1}v_{0}=0$, using (\ref{3.1}), we see that $b=a q$, then
\begin{align}
L_{-1}L_{1}(v_{j})=(p^{-j}[j]_{p,q}-ap^{-j}q^{j}-ap^{-j-1}q^{j+1})(p^{-j-1}[j+1]_{p,q}), j\in\Z,\label{3.23}\\
L_{1}L_{-1}(v_{j})=(p^{-j+1}[j-1]_{p,q}-ap^{-j+1}q^{j-1}-ap^{-j}q^{j})(p^{-j}[j]_{p,q}), j\in\Z,\label{3.24}
\end{align}
since the coefficient of $v_{j}$ in $p^{2j}q^{-2j}L_{-1}L_{1}$ (or $p^{2j}q^{-2j}L_{1}L_{-1}$) is a quadratic polynomial in $x:=q^{-j}[j]_{p,q}$ by Theorem 3.4, there is at most one other $v_{j}$ annihilated by $L_{-1}$. There are two subcases of $L_{-1}$ and $L_{1}$ as following:\\
subcase i: if dim(Ker$L_{-1})=2$, then Ker$(L_{-1})=\C v_{0}\oplus\C v_{-1}$;\\
subcase ii: either dim(Ker$L_{1})<2$, or dim(Ker$L_{-1})<2$.\\
\noindent{\it Proof of subcase i.~} If dim(Ker$L_{-1})=2$, we can assume $L_{-1}v_{-r}=0$ in addition to $L_{-1}v_{0}$, where $r$ is an integer. Putting $j=r, r-1, \cdots, 2, 1$ into
$$p^{-j}q^{j}L_{j}L_{-1}(v_{-r})-pq^{-1}L_{-1}L_{j}(v_{-r})=(\frac{[-1]_{p,q}}{p^{-1}}-\frac{[j]_{p,q}}{p^{j}})L_{j-1}(v_{-r}),\  j\in\Z,$$
we see that $L_{r-1}(v_{-r})=0, L_{r-2}(v_{-r})=0, \cdots, L_{1}(v_{-r})=0, L_{0}(v_{-r})=0$, then $p^{r}[-r]_{p,q}-ap^{r}q^{-r}=0$ by $L_{0}v_{k}=(p^{-k}[k]_{p,q}-ap^{-k}q^{k})v_{k}$, so $a=q^{r}[-r]_{p,q}$. It follows from (\ref{3.23}) that $(p^{r}[-r]_{p,q})(p^{r+1}[-r-1]_{p,q}-ap^{r+1}q^{-r-1}-ap^{r}q^{-r})=0$, then $r=1$ or $r=0$, hence Ker$(L_{-1})=\C v_{0}\oplus\C v_{-1}$.\\
\noindent{\it Proof of subcase ii.~}  Suppose that dim(Ker$L_{1})=2$, and dim(Ker$L_{-1})=2$, we have Ker$(L_{-1})=\C v_{0}\oplus\C v_{-1}$ and $r=1$ by proof of subcase i, then $a=-1, b=q$. If $L_{1}(v_{j})=0, j\in\Z$, by (6.24), $j$ satisfies $$(p^{-j-1}[j+1]_{p,q})(p^{-j}[j]_{p,q}-p^{-j}q^{j+1}[-1]_{p,q}-p^{-j-1}q^{j+2}[-1]_{p,q})=0,j\in\Z$$ then $j=-2$ or $j=-1$, hence Ker$(L_{1})=\C v_{-1}\oplus\C v_{-2}$. Next, we consider these two equations:
\begin{align}
pq^{-1}L_{-1}L_{2}(v_{j})-p^{-2}q^{2}L_{2}L_{-1}(v_{j})&=(\frac{[2]_{p,q}}{p^{2}}-\frac{[-1]_{p,q}}{p^{-1}})L_{1}(v_{j}), j\in\Z,\label{4.21}\\
p^{2}q^{-2}L_{-2}L_{1}(v_{j})-p^{-1}q^{1}L_{1}L_{-2}(v_{j})&=(\frac{[1]_{p,q}}{p^{1}}-\frac{[-2]_{p,q}}{p^{-2}})L_{-1}(v_{j}), j\in\Z,\label{4.22}
\end{align}
let $j=-1$ in (\ref{4.21}) and (\ref{4.22}), we have $L_{2}(v_{-1})=0$ and $L_{-2}(v_{-1})=0$, similarly, let let $j=0$ in (\ref{4.22}), we have $L_{-2}(v_{1})=0$, let $j=-2$ in (\ref{4.21}), we have $L_{2}(v_{-3})=0$. It follows that $\C v_{-1}$ is invariant under $V_{p,q}$ and hence is a submodule  of $V_{p,q}$, which contradicts the assumption that $V_{p,q}$ is indecomposable, hence this proves subcase ii.\\
Next, we consider the action of $L_{\pm2}$. Let $L_{-1}(v_{j})=m_{j}v_{j-1}$, where $m_{j}\in\C$, replacing $v_{j}$ by $$\frac{p^{-1}[1]_{p,q}p^{-2}[2]_{p,q}\cdots p^{-j}[j]_{p,q}}{m_{1}m_{2}\cdots m_{j}}v_{j},\ j\in\Z,$$ moreover, we have $m_{j}\neq0$ for $j\geq1$ , then
\begin{align}
L_{-1}(v_{j})=p^{-j}[j]_{p,q}v_{j-1},\  j\geq0,\label{4.23}
\end{align}
it follows from (\ref{3.23}) that
\begin{align}
L_{1}(v_{j})=(p^{-j}[j]_{p,q}-ap^{-j}q^{j}-ap^{-j-1}q^{j+1})v_{j+1},\  j\geq0.\label{4.24}
\end{align}
By (\ref{4.21}), (\ref{4.23}), (\ref{4.24}) and induction on $j$, we get
\begin{align}
L_{2}(v_{j})=(p^{-j}[j]_{p,q}-ap^{-j}q^{j}-ap^{-j-2}q^{j+1}[2]_{p,q})v_{j+2},\  j\geq0.\label{4.25}
\end{align}
Set
\begin{align}
L_{-2}(v_{j})=(f(j)+(p^{-j}[j]_{p,q}-ap^{-j}q^{j}-ap^{-j+2}q^{j+1}[-2]_{p,q}))v_{j-2},\ j\in\Z,\label{4.26}
\end{align}
using (\ref{4.22}), (\ref{4.23}), (\ref{4.24}) and (\ref{4.26}), we get
\begin{align}
p^{2}q^{-2}f(j+1)(p^{-j}[j]_{p,q}-ap^{-j}q^{j}-ap^{-j-1}q^{j+1})=
p^{-1}qf(j)(p^{-j+2}[j-2]_{p,q}-ap^{-j+2}q^{j-2}-ap^{-j+1}q^{j-1}), \label{4.27}
\end{align}
where $j\geq2$.
It is clear that there exists a $j_{0}\geq2$ such that
\begin{align}
p^{-j+2}[j-2]_{p,q}-ap^{-j+2}q^{j-2}-ap^{-j+1}q^{j-1}\neq0,\  j\geq j_{0}\geq2.\label{4.28}
\end{align}
We assume that $f(j)\neq0, j\geq j_{0}$, by (\ref{4.27}), we get
\begin{align}
f(j)=\frac{p^{-3(j-j_{0})}q^{3(j-j_{0})}(p^{-j_{0}+1}[j_{0}-1]_{p,q}-ap^{-j_{0}+1}q^{j_{0}-1}-ap^{-j_{0}}q^{j_{0}})}
{(p^{-j+1}[j-1]_{p,q}-ap^{-j+1}q^{j-1}-ap^{-j}q^{j})}\notag\\
\frac{p^{-j+2}[j-2]_{p,q}-ap^{-j+2}q^{j-2}-ap^{-j+1}q^{j-1}}{p^{-j+2}[j-2]_{p,q}-ap^{-j+2}q^{j-2}-ap^{-j+1}q^{j-1}},
 j\geq j_{0}.\label{4.29}
\end{align}
Using (\ref{4.25}),(\ref{4.26}) and (\ref{4.29}), we get
\begin{align}
p^{j}q^{-j}L_{-2}(v_{j})&=p^{j}q^{-j}(f(j)-(p^{-j}[j]_{p,q}-ap^{-j}q^{j}-ap^{-j+2}q^{j+1}[-2]_{p,q}))v_{j-2}\nonumber\\
&=\frac{F(x)}{G(x)}v_{j-2},\  j\geq j_{0},\\
p^{j}q^{-j}L_{2}(v_{j-2})&=p^{j}q^{-j}(p^{-j+2}[j-2]_{p,q}-ap^{-j+2}q^{j-2}-ap^{-j+4}q^{j-2}[-2]_{p,q})v_{j}\nonumber\\
&=(x+p^{2}[-2]_{p,q}-ap^{2}q^{-2}-ap^{4}q^{-2}[-2]_{p,q})v_{j},\  j\geq2,
\end{align}
where \\
\begin{align*}
F(x)&=(p^{-j_{0}+1}[j_{0}-1]_{p,q}-ap^{-j_{0}+1}q^{j_{0}-1}-ap^{-j_{0}}q^{j_{0}})
(p^{-j_{0}+2}[j_{0}-2]_{p,q}-ap^{-j_{0}+2}q^{j_{0}-2}-ap^{-j_{0}+1}q^{j_{0}-1})\\
&p^{3j_{0}}q^{-3j_{0}}f(j_{0})+(x-a-ap^{2}[-2]_{p,q})(x+p^{2}[-2]_{p,q}-ap^{2}q^{-2}-apq^{-1})
(x+p[-1]_{p,q}-apq^{-1}-a),\\
G(x)&=(x+p^{2}[-2]_{p,q}-ap^{2}q^{-2}-apq^{-1})
(x+p[-1]_{p,q}-apq^{-1}-a),
\end{align*}
and $x:=q^{-j}\{j\}$.

Since the coefficient of $v_{j}$ in $p^{2j}q^{-2j}L_{2}L_{-2}(v_{j})$ is a quadratic polynomial in $x$, this forces $f(j)=0, j\geq j_{0}$, hence
\begin{align*}
L_{-2}(v_{j})&=(p^{-j}[j]_{p,q}-ap^{-j}q^{j}-ap^{-j+2}q^{j+1}[-2]_{p,q})v_{j-2},\\
L_{-2}L_{2}(v_{j})&=(p^{-j}[j]_{p,q}-ap^{-j}q^{j}-ap^{-j-2}q^{j+1}[2]_{p,q})
(p^{-j-2}[j+2]_{p,q}-ap^{-j-2}q^{j+2}-ap^{-j}q^{j+3}[-2]_{p,q})v_{j},
\end{align*}
where $j\geq j_{0}$.
We have known that the coefficient of $v_{j}$ in $p^{2j}q^{-2j}L_{-2}L_{2}(v_{j})$ is a quadratic polynomial in $x$, and it has been identified as $(x-a-ap^{-2}[2]_{p,q})
(x+p^{-2}[2]_{p,q}-ap^{-2}q^{2}-aq^{2}[-2]_{p,q})$ from the proof of Theorem 3.4 for $j\in\Z$, so
\begin{align}
L_{-2}L_{2}(v_{j})&=(p^{-j}[j]_{p,q}-ap^{-j}q^{j}-ap^{-j-2}q^{j+1}[2]_{p,q})
(p^{-j-2}[j+2]_{p,q}-ap^{-j-2}q^{j+2}-ap^{-j}q^{j+3}[-2]_{p,q})v_{j},\label{4.32}
\end{align}
where $j\in\Z$.
\hfill$\Box$\vskip7pt
Using the subcases about $L_{-1}$ and $L_{1}$, we can derive the following two cases:
\begin{lemm}\label{2}
Either Ker$(d_{-1})=\C v_{0}$, or Ker$(d_{-1})=\C v_{0}\oplus\C v_{-1}$ and Ker$(d_{1})=0$.
\end{lemm}
\noindent{\it \textbf{Proof.}~}
Suppose Ker$(d_{-1})=\C v_{0}\oplus\C v_{-1}$ and Ker$(d_{1})=0$, then $\a=-1$ by proof of Lemma (6.6), so (\ref{4.32}) becomes
\begin{align}
L_{-2}L_{2}(v_{j})=&(p^{-j}[j]_{p,q}-ap^{-j}q^{j+1}[-1]_{p,q}-ap^{-j-2}q^{j+1}[-1]_{p,q})\notag\\
&(p^{-j-2}[j+2]_{p,q}-p^{-j-2}q^{j+3}-p^{-j}q^{j+3}[-2]_{p,q}[-1]_{p,q})v_{j},\  j\in\Z.\label{4.33}
\end{align}
 Let $j=-1$ in (\ref{4.21}), we have $L_{2}(v_{-1})\neq0$, let $j=-1$ and $j=-2$ in (\ref{4.33}), we have $L_{-2}(v_{1})=0$ and $L_{-2}(v_{0})\neq0$, it follows that (\ref{4.22}) fails for $j=0$, hence we get Ker$(L_{-1})=\C v_{0}$.
\hfill$\Box$\vskip7pt
It should be noted that all of our calculations are under the assumption that $L_{-1}(v_{0})=0$. By appropriately normalizing the $v_{j}^{'}s$ with $j$ negative, we can strength (\ref{4.23}) to
\begin{align}
L_{-1}(v_{j})=p^{-j}[j]_{p,q}v_{j-1}, j\in\Z,\label{4.34}
\end{align}
by (\ref{4.21}), we get
\begin{align}
L_{1}(v_{j})=(p^{-j}[j]_{p,q}-ap^{-j}q^{j}-ap^{-j-1}q^{j+1})v_{j+1}, \  j\neq{-1},\label{4.35}
\end{align}
by (\ref{4.21}), (\ref{4.25}), (\ref{4.34}) and (\ref{4.35}), and an induction on $j\leq-3$, we get
\begin{align}
L_{2}(v_{j})=(p^{-j}[j]_{p,q}-ap^{-j}q^{j}-ap^{-j-2}q^{j+1}[2]_{p,q})v_{j+2}, \  j\neq\{-1,-2\}.\label{4.36}
\end{align}
Let
\begin{align}
L_{1}(v_{-1})=Hv_{0},\ L_{2}(v_{-2})=Dv_{0},\ L_{2}(v_{-1})=Ev_{1},\label{4.37}
\end{align}
where $H, D, E\in\C$.
Put $j=-1$ in (\ref{4.21}) and using (\ref{4.34}) and (\ref{4.37}), we get
\begin{align}
q^{-1}E-p^{-1}q^{2}[-1]_{p,q}D=p^{-1}(q^{-1}[2]_{p,q}-q^{2}[-1]_{p,q})H.\label{4.38}
\end{align}
As for the equation $p^{-x}[x]_{p,q}-ap^{-x}q^{x}-ap^{-x-2}q^{x+1}[2]_{p,q}=0$, we make $j_{0}$ denote its solution when its integer solution exists, i.e.
\begin{align}
p^{-j_{0}}[j_{0}]_{p,q}-ap^{-j_{0}}q^{j_{0}}-ap^{-j_{0}-2}q^{j_{0}+1}[2]_{p,q}=0.\ \label{4.39}
\end{align}
Moreover, by (\ref{4.32}) and (\ref{4.36}), we can get
\begin{align}
L_{-2}(v_{j})=(p^{-j}[j]_{p,q}-ap^{-j}q^{j}-ap^{-j+2}q^{j+1}[-2]_{p,q})v_{j-2},\  j\neq\{1,0,j_{0}+2\}.\label{4.40}
\end{align}
Similarly, write
\begin{align}
L_{-2}(v_{0})=Fv_{-2},\   L_{-2}(v_{1})=Gv_{-1},\label{4.41}
\end{align}
put $j=-1$ and $j=-2$ in (\ref{4.32}), we can respectively get
\begin{align}
EG&=(p[-1]_{p,q}-apq^{-1}-ap^{-1}[2]_{p,q})
(p^{-1}-ap^{-1}q-apq^{2}[-2]_{p,q}),\label{4.42}\\
DF&=(p^{2}[-2]_{p,q}-ap^{2}q^{-2}-aq^{-1}[2]_{p,q})
(-a-ap^{2}q[-2]_{p,q}).\label{4.43}
\end{align}
Taking $j=-1$ and $j=1$ in (\ref{4.22}), and using (\ref{4.34}), (\ref{4.35}), (\ref{4.40}) and (\ref{4.41}), we have
\begin{align}
FH&=-pq^{-2}a(1+[2]_{p,q}a), \text{ when $j_{0}\neq{-3}\in\Z$,}\label{4.44}\\
GH&=-q^{-1}(1+[2]_{p,q}a)(a+p^{-1}), \text{ when $j_{0}\neq0\in\Z$.}\label{4.45}
\end{align}
Now we consider four cases.\\
Case 1: $a\neq0, -\frac{1}{p+q}, -p^{-1}$, then $j_{0}\neq-3,0$ by (\ref{4.39}), and $H\neq0$ by (\ref{4.44}). Replacing $v_{j}^{'}s$ with $j\leq-1$ by $(p[-1]_{p,q}-apq^{-1}-a)v_{j}$ and those $v_{j}^{'}s$ with $j\geq0$ by $Hv_{j}$, we get
\begin{align}
H=p[-1]_{p,q}-apq^{-1}-a.\label{4.46}
\end{align}
Then we have
\begin{align}
F=pq^{-1}a,\ G=a+p^{-1}\label{4.47}
\end{align}
from (\ref{4.44}), (\ref{4.45}) and (\ref{4.46}).
Applying (\ref{4.47}) to (\ref{4.42}) and (\ref{4.43}), we get
\begin{align}
D=p^{2}[-2]_{p,q}-ap^{2}q^{-2}-aq^{-1}[2]_{p,q},\ E=p[-1]_{p,q}-apq^{-1}-ap^{-1}[2]_{p,q}.\label{4.48}
\end{align}
In sum, we have
\begin{align}
L_{n}(v_{j})&=(p^{-j}[j]_{p,q}-ap^{-j}q^{j}-ap^{-j-n}q^{j+1}[n]_{p,q})v_{n+j},\  n=0,\pm1,2, j\in\Z.\label{4.49}\\
L_{-2}(v_{j})&=(p^{-j}[j]_{p,q}-ap^{-j}q^{j}-ap^{-j-n}q^{j+1}[-2]_{p,q})v_{j-2},\  j\neq j_{0}+2,\label{4.50}
\end{align}
we can check (\ref{4.50}) holds when $j=j_{0}+2$ by (\ref{4.32}). Hence, (\ref{4.49}) holds for all $n=0,\pm1,\pm2$ and all $j\in\Z$. Therefore, $V=V_{p,q}(a,b)$, where $b=aq$.\\
Case 2: $a=-\frac{1}{p+q}$, then $j_{0}\neq{-3}$ by (\ref{4.39}). By (\ref{4.38}), (\ref{4.42}), (\ref{4.43}) and (\ref{4.44}), we can deduce that
\begin{align}
H=0,\ q^{-1}E=p^{-1}q^{2}[-1]_{p,q}D,\ EG=\frac{p^{-2}q^{2}}{p^{2}-q^{2}},\ DF=\frac{p^{2}q^{-2}}{(p+q)^{2}}. \label{4.51}
\end{align}
Similarly, replacing $v_{j}^{'}s$ with $j\leq-1$ by $(p^{2}[-2]_{p,q}-ap^{2}q^{-2}-aq^{-1}[2]_{p,q})v_{j}$ and those $v_{j}^{'}s$ with $j\geq0$ by $Dv_{j}$, we get
\begin{align}
F=-\frac{pq^{-1}}{p+q},\ E=\frac{p^{-1}q}{p+q},\ G=\frac{p^{-1}q}{p-q}.\label{4.52}
\end{align}
Combining these with (\ref{4.34}), (\ref{4.35}), (\ref{4.36}) and (\ref{4.40}), we can get (\ref{4.49}) and (\ref{4.50}) also hold, similarly, we can check (\ref{4.50}) holds when $j=j_{0}+2$ by (\ref{4.32}). Hence, (\ref{4.49}) holds for all $n=0,\pm1,\pm2$ and all $j\in\Z$. Therefore, $V=V_{p,q}(a,b)$, where $b=-\frac{q}{p+q}, a=-\frac{1}{p+q}$.\\
Case 3: $a =-p^{-1}$, (\ref{4.39}) implies $j_{0}=-3$. By (\ref{4.42}) (\ref{4.43}) and (\ref{4.45}), we get
\begin{align}
DF=[-1]_{p,q},\ EG=0,\ GH=0,\label{4.53}
\end{align}
then by (\ref{4.38}), we have $G=0$.
Replacing $v_{j}^{'}s$ with $j\leq-1$ by $Fv_{j}$ and those $v_{j}^{'}s$ with $j\geq0$ by $(-a-ap^{2}q[-2]_{p,q})v_{j}$, we have $$F=-a-ap^{2}q[-2]_{p,q},$$ so $D=p^{-1}$. In conclusion, we get
\begin{align}
L_{n}(v_{j})=p^{-n-j-1}[n+j+1]_{p,q}v_{n+j},\ n=\pm1,\pm2,j\neq{-1}.\label{4.54}
\end{align}
In the previous hypothesis, and put $j=-1$ in equation (\ref{4.21}) and (\ref{4.22}) respectively, we have
\begin{align}
L_{1}(v_{-1})=Hv_{0},\ L_{2}(v_{-1})=p^{-2}q^{3}[-1]_{p,q}+p^{-2}[3]_{p,q}H,\ L_{-2}(v_{-1})=p^{3}[-3]_{p,q}+p^{3}q^{-3}H,\label{4.55}
\end{align}
Then we choose a appropriate complex number $\alpha$ such that $H=-q[-1]_{p,q}+[-1]_{p,q}[2]_{p,q}p^{-1}q\alpha$, so we get
\begin{align}
L_{n}(v_{-1})=(-q^{n}[-n]_{p,q}+[-n]_{p,q}[n+1]_{p,q}p^{-n}q^{n}\alpha)v_{n-1},\  n=\pm1,\pm2.
\end{align}
Then
$V_{p,q}(a,b)=V_{p,q}(\alpha)$, where
\begin{align}V_{p,q}(\alpha)=
\begin{cases}
L_{n}(v_{j})=p^{-n-j-1}[n+j+1]_{p,q}v_{n+j},& \text{$j\neq-1,$}\\
L_{n}(v_{-1})=(-q^{n}[-n]_{p,q}+[-n]_{p,q}[n+1]_{p,q}p^{-n}q^{n}\alpha)v_{n-1},& \text{$j=-1.$}
\end{cases}
\end{align}
Case 4: $a =0$, (\ref{4.39}) implies $j_{0}=0$. By (\ref{4.42}) (\ref{4.43}) and (\ref{4.44}), we get
\begin{align}
EG=[-1]_{p,q},\ DF=0,\ FH=0,\label{4.58}
\end{align}
then by (\ref{4.38}), we have $F=0$.
Replacing $v_{j}^{'}s$ with $j\leq-1$ by $Gv_{j}$ and those $v_{j}^{'}s$ with $j\geq0$ by $-p^{-1}$, we have $G=-p^{-1}$, so $E=p[-1]_{p,q}$. In conclusion, we get
\begin{align}
L_{n}(v_{j})=p^{-j}[j]_{p,q}v_{n+j},\ n,j=\pm1,\pm2,n+j\neq0.\label{4.59}
\end{align}
In the previous hypothesis, and put $j=-1$ and $j=1$ in equation (\ref{4.21}) and (\ref{4.22}) respectively, we have
\begin{align}
L_{1}(v_{-1})=Hv_{0},\ L_{2}(v_{-2})=p^{2}q^{-3}-p^{3}q[-3]_{p,q}H,\ L_{-2}(v_{2})=p^{-3}[3]_{p,q}+p^{-3}q^{3}H,\label{4.60}
\end{align}
Then we choose a appropriate complex number $\alpha^{'}$ such that $H=p[-1]_{p,q}+[-1]_{p,q}[2]_{p,q}pq^{-1}\alpha^{'}$, so we get
\begin{align}
L_{n}(v_{-n})=(p^{n}[-n]_{p,q}+p^{n}q^{-n}[-n]_{p,q}[n+1]_{p,q}\alpha^{'})v_{0},\  n=\pm1,\pm2.\label{4.61}
\end{align}
Then
$V_{p,q}(a,b)=V_{p,q}^{'}(\alpha^{'})$, where
\begin{align}V_{p,q}^{'}(\alpha^{'})=
\begin{cases}
L_{n}(v_{j})=p^{-j}[j]_{p,q}v_{n+j},& \text{$j\neq-n,$}\\
L_{n}(v_{-n})=(p^{n}[-n]_{p,q}+p^{n}q^{-n}[-n]_{p,q}[n+1]_{p,q}\alpha^{'})v_{0},& \text{$j=-n.$}\label{4.62}
\end{cases}
\end{align}
Then, when $L_{-1}$ annihilates some $v_{j},j\in\Z$, we can get three $V_{p,q}$-modules:
$V_{p,q}(a,b)$, $V_{p,q}(\alpha)$, $V_{p,q}^{'}(\alpha^{'})$.\\
\textbf{Case II:}\\
Similarly to the discussion of $L_{-1}$ annihilates, we can get other modules of $V_{p,q}$, that is,
\begin{align}V_{p,q}(\beta)=
\begin{cases}
L_{n}(v_{j})=-q^{n+j-1}[-n-j+1]_{p,q}v_{n+j},& \text{$j\neq1,$}\\
L_{n}(v_{1})=(-q^{n}[-n]_{p,q}+p^{-n}q^{n}[n]_{p,q}[-n+1]_{p,q}\alpha)v_{n+1},& \text{$j=1,$}
\end{cases}
\end{align}
 and
\begin{align}V_{p,q}^{'}(\beta^{'})=
\begin{cases}
L_{n}(v_{j})=p^{-j}[j]_{p,q}v_{n+j},& \text{$j\neq-n,$}\\
L_{n}(v_{-n})=(p^{n}[-n]_{p,q}+p^{n}q^{-n}[n]_{p,q}[n+1]_{p,q}\alpha^{'})v_{0},& \text{$j=-n.$}
\end{cases}
\end{align}
Now we have the following theorem:
\begin{theo}
Let $q$ be not a root of unity, if $V$ is an indecomposable Harish-Chandra $V_{p,q}$-module with one dimensional weight spaces, then $V$ is isomorphism to one of the following $V_{p,q}$-module:
\begin{align*}V_{p,q}(a,b): &\ L_{n}(v_{j})=(p^{-k}[k]_{p,q}-ap^{-k}q^{k}-bp^{-k-j}q^{k}[j]_{p,q})v_{n+j},\\
V_{p,q}(\alpha):
&\begin{cases}
L_{n}(v_{j})=p^{-n-j-1}[n+j+1]_{p,q}v_{n+j},& \text{$j\neq-1,$}\\
L_{n}(v_{-1})=(-q^{n}[-n]_{p,q}+[-n]_{p,q}[n+1]_{p,q}p^{-n}q^{n}\alpha)v_{n-1},& \text{$j=-1.$}
\end{cases}\\
V_{p,q}^{'}(\alpha^{'}):
&\begin{cases}
L_{n}(v_{j})=p^{-j}[j]_{p,q}v_{n+j},& \text{$j\neq-n,$}\\
L_{n}(v_{-n})=(p^{n}[-n]_{p,q}+p^{n}q^{-n}[-n]_{p,q}[n+1]_{p,q}\alpha^{'})v_{0},& \text{$j=-n.$}\label{4.62}
\end{cases}\\
V_{p,q}(\beta):
&\begin{cases}
L_{n}(v_{j})=-q^{n+j-1}[-n-j+1]_{p,q}v_{n+j},& \text{$j\neq1,$}\\
L_{n}(v_{1})=(-q^{n}[-n]_{p,q}+p^{-n}q^{n}[n]_{p,q}[-n+1]_{p,q}\alpha)v_{n+1},& \text{$j=1.$}
\end{cases}\\
V_{p,q}^{'}(\beta^{'}):
&\begin{cases}
L_{n}(v_{j})=p^{-j}[j]_{p,q}v_{n+j},& \text{$j\neq-n,$}\\
L_{n}(v_{-n})=(p^{n}[-n]_{p,q}+p^{n}q^{-n}[n]_{p,q}[n+1]_{p,q}\alpha^{'})v_{0},& \text{$j=-n.$}
\end{cases}
\end{align*}
where $\a,\b,\alpha^{'},\beta^{'},a,b\in\C$.
\end{theo}

\end{document}